\documentclass{article}

\usepackage{theorem}

\makeatletter
\@addtoreset{equation}{subsection}
\@addtoreset{equation}{section}
\makeatother
\newtheorem{theorem}[equation]{Theorem}
\newtheorem{proposition}[equation]{Proposition}

{\theorembodyfont{\rmfamily}

}

\newcommand{\C}{{\mathbf C}}
\newcommand{\Q}{{\mathbf Q}}
\newcommand{\Z}{{\mathbf Z}}
\newcommand{\R}{{\mathbf R}}
\newcommand{\p}{\mathbf P}
\newcommand{\pn}{\mathbf P^n}
\newcommand{\Ol}{{\cal O}}
\newcommand{\Li}{{\cal L}}
\newcommand{\opn}{{\Ol_{\mathbf P^n}}}
\newcommand{\f}{\varphi}
\newcommand{\ra}{\longrightarrow}

\begin{document}

\title{Characterization theorems for the projective space
and vector bundle adjunction.}

\author{Marco Andreatta}

\date{}
\maketitle

\begin{abstract}
We consider some conditions under which a smooth projective variety $X$ 
is actually the projective space. 
We also extend to the case of positive characteristic some
results in the theory of vector bundle adjunction.
We use methods and techniques
of the so called Mori theory, in particular the study of rational
curves on projective manifolds.
\end{abstract}


\section{Introduction}
\label{intro}
Let $X$ be a smooth projective variety of dimension $n$, defined over
an algebraically closed field $k$; we denote by
$TX$ its tangent bundle and by $K_X = - det( TX)$ its canonical bundle.

A natural problem is to find simple conditions under which 
the manifold $X$ is actually the projective space $\pn$.
A very famous one is given by the following 
theorem of S. Mori (see \cite{RefMo79}).

\noindent
{\bf Theorem}. 
{\it $X$ is the projective space if and only if $TX$ is ample.}

In the last year it has been generalized, over the
field of complex number, in the following two directions.

\noindent
{\bf Theorem} (see \cite{RefCMS01} and \cite{RefKe01}).
{\it ($k = \C$). $X$ is the projective space 
if and only if 
$-{K_X }^.C \geq (n+1) $ 
for any curve $C \subset X$.}

\noindent
{\bf Theorem} (see \cite{RefAW01}). 
{\it ($k = \C$).  $X$ is the projective space if and only if
there exists an ample subsheaf $E \subset TX$. Of course we then have that
$E = T \pn$ or $E = \oplus^r \opn (1)$.}

The actual proofs of the last two theorems are valid only in the case $k = \C$,
although they may be valid even in positive characteristic.
In the positive characteristic case we have the following.

\noindent
{\bf Theorem} (see \cite{RefKK00}). 
{\it Let $H$ be an ample line bundle on $X$. Then
either $K_X +(n-1)H$ is nef or the pair $(X,H)$ is one of the following:
$X = \pn$ and $H$ is a hyperplane section, $X$ is a quadric $\Q \subset \p ^{(n+1)}$
and $H$ is a hyperplane section, $X$ is the projectivization of a rank $n$ vector bundle over 
a smooth curve $A$ and $H$ is the tautological bundle.}

The aim of this note is to give some other simple characterizations of projective space
(or even of quadrics and scrolls) which are slight generalizations of the above theorems.

First of all we have the following more general form of
the theorem in \cite{RefAW01} 
(to obtain it consider $\f$ to be the constant map).

\begin{theorem} 
\label{theoA}  ($k = \C$). Let  $\f : X \ra Y$ be a surjective map to a normal 
projective variety $Y$ of lower dimension (i.e. $m:= dim Y < n$) and let 
$T_{X/Y}$  denote the relative tangent bundle.
\item{i)} If $T_{X/Y}$ contains an ample vector bundle $E$
then $X = \pn$, $dim Y = 0$ and $E = T\pn$ or $E = \oplus^r \opn (1)$. 
\item{ii)} If $\f$ is equidimensional and smooth and $T_{X/Y}$ contains a $\f$-ample vector bundle $E$
then $\f : X \ra Y$ is a scroll, i.e. $X$ is the projectivization of a rank $n-m +1$ vector bundle over 
$Y$.
\end{theorem}

Part i) of the above theorem solves the 
problem 4.4 in \cite{RefPe01}.

As for the theorem proved in \cite{RefKK00}, we can state a more general version
using a vector bundle instead of a line bundle (to get the above quoted version
one simply takes $E = \oplus^r H$). 

\begin{theorem}
\label{theoB}
Let $E$ be an ample vector bundle on $X$ of rank $r$.
\item{i)} If $r \geq n$ then $K_X + detE$ is nef unless $ X = \pn$ and 
$E = \oplus^n \opn (1)$ 
\item{ii)} If $r = n-1$ then either $K_X + det E$ is nef or 
the pair $(X,E)$ is one of the following:
\item{a)} $X = \pn$ and $E = \opn (2)\oplus (\oplus^{n-2} \opn (1))$ ,
\item{b)} $X = \pn$ and $E = \oplus^{n-1} \opn (1) $,
\item{c)} $X$ is a quadric $\Q \subset \p ^{(n+1)}$
and  $E = \oplus^{n-1} \Ol(1)$, 
\item{d )} $X$ is the projectivization of a rank $n$ vector
bundle over  a smooth curve $A$, $\f : X \ra A$, and $E= \f^*(\f_*(E(-1)))\otimes \Ol(1)$.
\end{theorem}

Under the additional hypothesis
that $E$ is spanned by global sections (part of) the theorem 
has been proved in \cite{RefSa00}; it has been proved 
in the case $k = \C$ in \cite{RefYZ90}.

In the last section we state a natural
conjecture, which is true if $k = \C$, and, as a test, we prove an easy
consequence of it.

I would like to thank the referee for pointing out many inaccuracies in the 
first version of the paper.

\section{Proofs}
\label{sec:1}
In this section we prove the theorems stated in the introduction. Our notation is 
consistent with the standard in use in algebraic geometry in particular with the one used 
in the book \cite{RefKo96} to which we frequently refer the reader.

\noindent 
{\bf Proof} (of Theorem \ref{theoA}, part i). 
We have to prove that $dim Y = 0$, then the theorem 
follows from the theorem in \cite{RefAW01}; 
let us then assume by contradiction that
$m:= dim Y > 0$. 

Take the scheme theoretic intersection of $m-1$ general sections of a line bundle on
$X$ which is the pull back of a very ample line bundle on $Y$; call it $X'$
and call $\f'$ and $Y'$ the restriction
of $\f$ to $X'$ and its image respectively.
Note that $dim Y' = 1$ and that, by Bertini' s theorem, 
$X'$ and $Y'$ are smooth. 
Thus  $\f' : X' \ra Y'$ is  equidimensional and 
moreover $\f'$ is a map as in the theorem,
namely $T_{X'/Y'}$ contains an ample vector bundle $E' :=E_{X'}$ 
Let $F$ be the general fiber 
of this map.
By the fact that ${T_{X'/Y'}}_{|F} = T_F$ and by the 
theorem in \cite{RefAW01},  $F$ is isomorphic to $\p^{(n'-1)}$,
$n' = dim X'$. 
Thus, by Fuijta's characterization of scrolls, 
\cite{RefFu87} lemma (2.12), $\f'$ is a $\p^{(n'-1)}$-bundle. This is a contradiction with 
Lemma (1.2) of \cite{RefCP98} 
(which says exactly that there is no ample subbundle of the relative tangent
bundle of a $\p$-bundle).

\noindent 
{\bf Proof} (of Theorem \ref{theoA}, part ii). 
The proof is by induction on $m: =dim Y$. 
If $m = 0$ then this is the theorem in \cite{RefAW01}. 
Let us assume that the theorem is true for $m-1$. 
Then the locus of points in $Y$ over
which the fiber  is not $\p^{r}$, call it $D$, is discrete. 
In fact, take a general very ample divisor $A$ of $Y$ and let 
$X' = \f^{-1} (A) \ra A $; this map is by induction a $\p^r$-bundle
and therefore $D \subset Y\setminus A$.
Thus we can apply Lemma (3.3) of \cite{RefAM97} 
to have an ample line bundle on $X$ 
whose restriction to a general fiber is $\Ol(1)$. 
The theorem follows from the
Fuijta characterization of scrolls, 
\cite{RefFu87} lemma (2.12). 

\noindent 
{\bf Proof} (of Theorem \ref{theoB}, part i)).  Let $r \geq n$ and assume that 
$K_X +det E$ is not nef. Then, by the Mori cone theorem, there exists a rational curve
$f :\p^1 \ra C\subset X$ such that 
$$(K_X +det E)^. C < 0, \ \ \  {-K_X} ^.C \leq (n+1);$$
moreover for every rational curve $C$, by the ampleness of $E$,
we have that 
$detE^. C \geq rank(E)(= r \geq n)$.
These inequalities actually imply  
$${-K_X} ^.C = (n+1), \ \ \ \ detE^. C = r  = n.$$
Let $V\subset Hom(\p^1, X)$ be a closed irreducible component containing
$f$ (see \cite{RefKo96} for this and the subsequent definitions).
Since $detE^. C = r$ and $E$ is ample, 
$V$ is actually an unsplit family of rational curve 
(see  \cite{RefKo96}, definition IV  2.1). 

Let $0\in \p^1$ be a point, $x =f(0)$ and 
$V_x := V \cap Hom(\p^1, X,0 \mapsto x)$. 
Then $dim_{f} V_x \geq {-K_X} ^.C  = (n+1)$
(see \cite{RefKo96}, II.1) and
therefore the locus of the curve in $V_x$
covers all of $X$, i.e.  $Locus(V_x) = X$
(note that to go from homomorphisms to curves,
i.e. from $Hom$ to $Chow$, we have to quotient out by 
$Aut(\p^1, 0)$ which has dimension $2$).

In particular this means that $X$ is 
{\it rationally chain connected} with respect to
the family $V$, 
i.e. for every two points on $X$ there is 
a chain of rational curves parameterized by morphisms from $V$
which joins these two points
(we have actually proved that $X$ is 
{\it rationally connected} with respect to
the family $V$, 
i.e. for every two points on $X$ there is 
a  rational curves parameterized by morphisms from $V$
which joins them). 

These varieties have many nice properties, here and 
later on we will use the one described in 
the Proposition IV.3.13.3 of \cite{RefKo96}.
Namely it says that $A_1(X)_{\Q}= A_1(X)_{\Z} \otimes _{\Z} \Q$ is generated by the
irreducible components of fibers of the universal family of $V$
(we refer to the section II.4.1 of
\cite{RefKo96} for the precise definition of $A_1(X)_{\Z}$,
respectively of $B_1(X)_{\Z}$ or of $N_1(X)_{\Z}$, 
the quotient of the group of
1-cycles by the subgroup of 1-cycles rationally, respectively algebraically
or numerically, equivalent
to the  zero cycle). 

Since $V$ is an unsplit family this implies that the fibers
of the universal family of $V$ are irreducible and moreover that they are
all algebraically equivalent; therefore $B_1(X)_{\Q} = \Q$.
Since $N_1(X)_{\Z}$ is a quotient of
$B_1(X)_{\Z}$, we have that $\rho (X) =: dim N_1(X)_{\R} = 1$.
 
Let $H := -(K_X +det E)$.  Since $H^.C =1$, $H$ is an ample line bundle
and $K_X + nH$ is not nef. Thus we can apply the theorem in \cite{RefKK00} 
to the pair $(X,H)$ to get 
$X = \p^n$. The description of $E$ follows from the characterization
of uniform vector bundles on $\p^n$, see \cite{RefEi80}.

\noindent 
{\bf Proof} (of Theorem \ref{theoB}, part ii))  Let $r = n-1$ and assume that 
$K_X +det E$ is not nef. Then, by the Mori cone theorem, there exists a rational curve
$f :\p^1 \ra C\subset X$ such that 
$$(K_X +det E)^. C < 0, \ \ \  {-K_X} ^.C \leq (n+1), \ \ \ \ detE^. C \geq r = n-1.$$
These inequalities imply the following two possibilities: 

\noindent
$1)\ \ \ {-K_X} ^.C = n, \ \ \ \ detE^. C = n-1,\ \ \ \ \ \ \ \ \ \ \ \ $

\noindent
$2)\ \ \ {-K_X} ^.C = (n+1), \ \ \ \ detE^. C  = n\ \ \hbox{or}\ \  (n-1)$

Let $V\subset Hom(\p^1, X)$ be a closed irreducible component containing
$f$ and let $H := -(K_X +det E)$. Note that in both cases
the family $V$ is unsplit.

In the first case we will adapt the proof of \cite{RefKK00}. 
Namely through every $x\in X$ we have an unsplit family of rational
curves, associated to $V_x$, of dimension at least $(n-2)$. If there is an $x$ such
that these curves cover all $X$, then as in the proof of part i), we obtain that $X = \p^n$
and thus we get to a contradiction. 

Thus we can assume that for every $x$ the curves 
$\{f_t(\p^1): f_t \in V_x \}$ sweep out a divisor $B_x$.
These divisors form an algebraic family for $x$ in a suitable open
set $X^0 \subset X$.

If for any $x_1, x_2 \in X$ we have $B_{x_1}\cap B_{x_2} \not= \emptyset$, 
then any two points of $X$ are connected
by a chain of length $2$ of the form $C_{t_1}\cap C_{t_2}$, with $C_t$ 
in the family $V$.
We can apply again the Proposition  IV.3.13.3 of \cite{RefKo96}
and we have $\rho (X) = 1$. 
Therefore $H$ is ample and  
$K_X + (n-1)H$ is not nef; thus we can apply the theorem in \cite{RefKK00} 
to the pair $(X,H)$ to get that $X$ is a quadric.
The description of $E$ follows easily from the fact that $E_l = \oplus^{n-1} \Ol(1)$
for all lines on the quadric.

If for general $x_1, x_2 \in X$ we have $B_{x_1}\cap B_{x_2} = \emptyset$ 
then one can apply 
exactly the same argument of the last point (point 13) of the proof in \cite{RefKK00}. 
Namely it can be seen that the algebraic family of divisors $B_x$ determines
a morphism $p: X \ra A$ to a smooth curve $A$ and that this is actually a 
$\p^{(n-1)}$-bundle. 
Note that we will use the line bundle $H$ which is a priori not ample in our case,
but it is ample on a general fiber $Y$ of the map $p$ since $\rho (Y) = 1$ 
(again by the Proposition  IV.3.13.3 of \cite{RefKo96},
since $Y$ is rationally connected with respect to the family
$V$). 
The given description of $E$ follows easily since
$E_Y = \oplus^{n-1} \Ol(1)$. 

Let us then consider the second case. For every $x\in X$ the family $V_x$
is of dimension $(n+1)$. Then the locus of the curve in $V_x$
covers all of $X$ and as above, by the
the Proposition  IV.3.13.3 of \cite{RefKo96}, we have that
$\rho (X) = 1$. 

If  $detE^. C  = n$ then $H. C  = 1$; therefore
$H$ is ample and $K_X + (n+1)H$ is not nef.
The theorem in \cite{RefKK00}  applies  (or apply the  
part i)) and we obtain that $X = \p^n$. Moreover
$E = \Ol(2)\oplus (\oplus^{n-2} \Ol(1))$ (see \cite{RefEi80}).

If  $detE^. C  = (n- 1)$ then $H. C  = 2$, thus we cannot apply the results in 
\cite{RefKK00}. The rest of the section deals with this case.

First let us
consider the projectivisation $p: \p(E)\ra X$ with the relative $\Ol(1)$
bundle which we will denote by $\Li$. Note that since
$detE^. C  = (n- 1)$ for 
$f :\p^1 \ra C\subset X$ we have $f^*E=\oplus^{n-1}\Ol(1)$;
thus $\p(f^*E)=\p^1\times\p^{n-2}$.
It was observed in \cite{RefAW01} that in this situation,
for any $f\in V$ and $y\in p^{-1}(f(0))$, 
we have a unique lift-up $\hat f: \p^1\ra \p(E)$
such that $p\circ\hat f=f$, $deg(\hat f^*(\Li))=1$ and $\hat f(0)=y$.
Namely the morphism $\hat f$ is
obtained by simply composing  $\p(f^*E)\ra\p(E)$ with the product morphism
$\p^1\ra\p^1\times\{y\}\subset\p^1\times\p^{n-2} $.
Moreover, for a generic $f$ we have 
$\hat f^*T\p(E)=f^*TX\oplus\Ol^{\oplus (n-2)}$. 

We can choose an irreducible $\hat V\subset Hom(\p^1,\p(E))$
which parameterizes these lift-ups, that is, via the natural morphism
$p_*: Hom(\p^1,\p(E))\ra Hom(\p^1,X)$ defined by $p_*(\hat f)=p\circ
\hat f$, the component $\hat V$ dominates $V$.

\begin{proposition} (see \cite{RefAW01}) The morphism $p_*: \hat V\ra V$ is proper and
thus surjective; moreover $\hat V$ is an unsplit family.
\end{proposition}

\noindent 
{\bf Proof.} The proof in \cite{RefAW01} is in fact characteristic free.

The family $\hat V$ defines a relation of {\sl rational chain connectedness with
respect to} $\hat V$, which we shall call rc$\hat V$ relation for short, in the
following way: $x_1,\ x_2\in \p(E)$ are in the rc$\hat V$ relation if there
exists a chain of rational curves parameterized by morphisms from $\hat V$
which joins $x_1$ and $x_2$. We use the following result of Campana and
Koll\'ar-Miyaoka-Mori.

\begin{theorem} (see \cite{RefKo96}, IV.4.16). There exists an open subset
$\p^0\subset \p(E)$ and a proper surjective morphisms with connected fibers
$\f^0: \p^0\ra Z^0$ onto a normal variety, such that the fibers of $\f^0$
are equivalence classes of rc$\hat V$ relation. 
We shall call the morphism $\f^0$ an rc$\hat V$ fibration.
\end{theorem}

Note that a general fiber of the fibration $\f^0$, 
$Y\subset \p(E)$, 
is rationally chain 
connected with respect to the
family $\hat V$. In particular, as above by the Proposition IV.3.13.3
in \cite{RefKo96}, $\rho(Y)=1$.

By the surjectivity of $p_*:
\hat V\ra V$ and rational connectedness of $X$, the restriction map
$p_Y: Y\ra X$ is surjective. Since $\rho(Y)=1$, it has no positive
dimensional fiber, so it is a finite morphism. 

$Y$ is projective variety of dimension $n$ which is a locally complete intersection
such that $N_{Y/\p(E)} = \Ol_Y$.  In particular we have that 
the canonical sheaf of $Y$ is defined and it is actually
$K_Y = ({K_{\p(E)}})_{|Y}$.

Since $\Li ^.\tilde C = 1$ for any curve $\tilde C$ from $\hat V\cap Hom(\p^1,Y)$
and since $K_{\p(E)} = p^* K_X + p^*(det E)- (n-1)\Li$,
modulo numerical equivalence we have 
$-K_Y = -({K_{\p(E)}})_{|Y} \equiv (n+1)\Li$.

Let $y \in Y$ be a general point 
and let $\hat C $ be a general curve in the family such that if 
$f: \p^1 \ra Y \subset \p(E)$ is the normalization of $\hat C$ then $y = f(0)$
for $0 \in \p^1$. 

$\hat V_y =  \hbox{Hom}(\p^1,\p(E), 0 \ra y) \cap \hat V$ is
contained in $Y$ and it has dimension 
$dim \hat V_y  \geq -K_{\p(E)}.\hat C = (n+1)$.
Therefore there exists an $(n-1)$-dimensional family of
rational curves through $y$, corresponding to $V_y$.

Let $\pi: \tilde Y \ra Y $ be the normalization of $Y$, $\tilde \Li$ the pull back of
$\Li$ and $\tilde V_y$ the pull back of the family $\hat V_y$. This is again
an $(n-1)$-dimensional family of rational curves through $y\in \tilde Y $
and $\pi^* \Li .\tilde C = 1$ for any curve in the family

Note that $-K_{\tilde Y} = \pi^* (-K_Y) + ( \hbox{conductor of } \pi)$ where the 
conductor of $\pi$ is an effective divisor which is zero
iff $\pi$ is an isomorphism. Thus 
$$ - (K_{\tilde Y}. (\pi^* \Li)^{n-1})\geq  - (K_{Y}.  \Li^{n-1}) = (n+1)\Li^{n} = 
(n+1)(\pi^* \Li)^{n}, $$
the first inequality coming from the fact that the conductor is effective
and by projection formula.

These data fulfill 
the assumption of Lemma 8 in \cite{RefKK00} which implies that  
$\tilde Y = \pn$ and $\pi^* \Li = \Ol(1)$.

But, since 
$-K_{\tilde Y} = -(n+1)(\pi^* \Li) + ( \hbox{conductor of } \pi)$, 
this implies that the conductor
is zero and thus that $ Y = \tilde Y = \pn$ and $\Li = \Ol(1)$.

To conclude we consider the finite morphism
$p_Y: Y  = \pn \ra X$; we claim that it is separable and thus that 
it is an isomorphism by Lazarsfeld theorem
(see \cite{RefKo96}, V.3.5).

To prove this it is enough to see that there is a point $y\in Y$ 
for which no non zero vector in $T_yY$ goes to zero
via the differential $(d p_Y)_y$.
Let then $x \in X$ be a general point and let $y = p_Y^{-1}(x)$.

We recall that we are dealing 
with an unsplit family $V$ such that for every
$f \in V$ we have $deg f^*H =2$, where $H$ is the ample line bundle defined
above.
Thus we can apply the 
theorem 3.3 in \cite{RefKe02} which says that if
$x$ is a general point of $X$,  for any $f \in V_x$ its image
is smooth at $x$.

On the other hand, since $Y  = \pn$,  for every $0 \not= v \in T_yY$ 
we have a   $\hat f \in \hat V_y$ with $(d \hat f)_0({\partial \over \partial t}) = v$ 
and also
$f = p_Y\circ \hat f$. This implies easily our claim.

\section{Comments and remarks}
\label{sec:2} 
We want to make some comments on the above theorems. 

First of all we think that the assumption of smoothness of the map 
$\f$ in the theorem \ref{theoA} ii) is probably redundant and the theorem should be
valid under the weak assumption that 
$\f$ is simply equidimensional and surjective into a 
normal variety. 
The problem is to show the smoothness of $Y$ under 
these weaker hypothesis. 

\smallskip
The last part of the proof of Theorem \ref{theoB} is based on the proof,
given in \cite{RefAW01}, of the following.

\begin{proposition} (see \cite{RefAW01} proposition (1.2)) ($k = \C$)
\label{AW}
Let $X$ be a smooth projective manifold, $E$
a vector bundle of rank $r$ on $X$ and let 
$V\subset Hom(\p^1, X)$ be a closed irreducible component
which is unsplit. Assume that $X$ is rationally chain connected with
respect to $V$ and that for any $f\in V$ we have
$f^* E = \oplus^r \Ol (a)$. 
Then there exists a (uniquely defined) line bundle $L$ over
$X$ such that $f^* L = \Ol (a)$ and $ E = \oplus^r L$.
\end{proposition}

\noindent
{\bf Conjecture}. 
The proposition \ref{AW} is true over any algebraic
closed field $k$.

By the proposition many results on adjunction theory of vector bundles 
would follow easily from the, usually known, case $r= 1$. 
As an example let us give the following generalization
of Theorem \ref{theoB}. i) 
(it has been proved under the condition $r \leq (n+1)$
in \cite{RefOh99}).

\begin{theorem} Let $E$ be an ample vector bundle on $X$ of rank $r$,
let $$\tau (E) := \mathrm{ min}\{t \in \R: K_X + t (det E) \mathrm{\ \ is \ nef\ }\}$$ be
its nef value.  Then $\tau \leq {n+1\over r}$ with equality if and only if $ X = \pn$ and 
$E = \oplus^r \Ol(1)$.
\end{theorem}

\noindent 
{\bf Proof.}
Let $\tau$ be as in the theorem and let
$f :\p^1 \ra C\subset X$  be a rational curve such that 
$$(K_X + \tau det E)^. C = 0, \ \ \  {-K_X} ^.C \leq (n+1), \ \ \ \ detE^. C \geq r$$
(it exists by the Mori cone theorem).

These inequalities imply that $\tau \leq {n+1\over r}$, with equality if and only if
${-K_X} ^.C = (n+1)$ and $detE^. C = r$
(at this point if $k =\C$ we could use the theorem in  \cite{RefCMS01}
or in \cite{RefKe01}) .
The above proposition thus applies (with $V$ a component of 
the $Hom$ scheme which contains 
$f :\p^1 \ra C$) and it gives that $E = \oplus^r L$,
for an ample line bundle $L$.
Thus $-K_X = (n+1)L$ and the theorem follows
from the theorem \ref{theoB}.

Similarly one can easily describe the pair $(X, E)$ with $E$ an ample vector bundle of rank 
$r=n$ for which $K_X +det E$ is (nef but) not ample 
(if $k = \C$ this was done in \cite{RefFu90}).

\medskip
{\small 
\noindent
Marco Andreatta, 
Dipartimento di Matematica, 
Universit\`a di Trento ,
I-38100 Povo (TN), Italy.  

\noindent
e-mail: andreatt@science.unitn.it

\end{document}
